\documentclass{amsart}

\usepackage{color}
\usepackage{amsthm,amssymb,verbatim}
\usepackage{graphicx}
\usepackage{enumerate}

 \newcommand{\Dd}{\mathcal{D}}

 \newcommand{\RR}{\mathbf{R}}  % reals
 \newcommand{\BB}{\mathbf{B}}  %ball
 \newcommand{\Tan}{\operatorname{Tan}}

\input epsf
\def\begfig {
\begin{figure}
\small }
\def\endfig {
\normalsize
\end{figure}
}

    \newtheorem*{claim*}{Claim}
    \newtheorem*{theorem*}{Theorem}

    \newtheorem{theorem}    {Theorem}     % [section]
    
    \newtheorem{corollary}  [theorem]     {Corollary}
    \newtheorem{proposition}       [theorem]       {Proposition}

    \theoremstyle{definition}

    \theoremstyle{definition}
    \newtheorem{remark}   [theorem]       {Remark}
    \newtheorem*{remark*}{Remark}

\usepackage{hyperref}
\usepackage[alphabetic, msc-links]{amsrefs}

\begin{document}

\setlength{\baselineskip}{1.2\baselineskip}

\renewcommand{\thesubsection}{\thetheorem}
   % this is so subsections and theorems, etc will be
   % numbered together.

\title[Curvatures of Embedded Minimal Disks]{Curvatures of embedded minimal disks
   blow up on subsets of $C^1$ curves}
\author{Brian White}
\address{Mathematics Department, Stanford University, Stanford, CA 94305}
\thanks{The author would like to thank David Hoffman for helpful suggestions.
The research was partially supported by NSF
  grant~DMS-1105330}
\email{white@math.stanford.edu}
\date{March 28, 2011.  Revised February 19, 2015.}
\begin{abstract}
Using results of Colding-Minicozzi and an extension due to Meeks, 
we prove that a sequence of properly embedded minimal disks in a $3$-ball
must have a subsequence whose curvature blow-up set lies in a union of disjoint $C^1$ curves.
\end{abstract}


\begin{comment}
\end{comment}

\maketitle

\section*{Introduction}\label{section:intro}

\newcommand{\LL}{\mathcal{L}}     %lamination
\newcommand{\HH}{\operatorname{\mathcal H}}
\newcommand{\HJ}{\operatorname{\mathcal H}_J}
\newcommand{\DD}{\operatorname{\bf\mathcal D}}   % horizontal disk
\newcommand{\CylJ}{\BB(0,1)\times z(J)}                    %cylinder

Let $D_n$ be a sequence of minimal disks that are properly embedded in an open
subset $U$ of $\RR^3$ or more generally of a $3$-dimensional Riemannian manifold. 
By passing to a subsequence, we may assume that there is a relatively closed
subset $K$ of $U$ such that the curvatures of the $D_n$ blow up at each point of $K$
(i.e., such that for each $p\in K$, there are points $p_n\in D_n$ converging to $p$ such that curvature of $D_n$
at $p_n$ tends to infinity as $n\to\infty$) and such that $D_n\setminus K$ converges
smoothly on compact subsets of $U\setminus K$ to a minimal lamination $L$ of $U\setminus K$.
It is natural to ask what kinds of singular sets $K$ and laminations $L$ can arise in this way.
In this paper, we prove:

\begin{theorem}\label{MainFromIntro} Every point of $K$ contains a neighborhood $W$ such that 
 $K\cap W$ is (after a rotation of $\RR^3$) contained in the graph of a $C^1$ function from $\RR$
 to $\RR^2$.  
\end{theorem}
 
This extends previous results of Colding-Minicozzi and of Meeks.
In particular, if one replaces ``$C^1$'' by ``Lipschitz" in Theorem~\ref{MainFromIntro}, then the result is implicit
in the work of Colding and Minicozzi.
(See~\cite{ColdingMinicozziIV}*{Section I.1} and~\cite{ColdingMinicozziIV}*{Theorem 0.1} for a very
similar result.)
  Thus if $K$ is a curve,
it must be a Lipschitz curve.
Meeks later showed that if $K$ is a Lipschitz curve,
then it must be a $C^{1,1}$ curve~\cite{MeeksRegularity}.

\begin{comment}
Colding and Minicozzi also showed that if $U=\RR^3$ (the ``global case'') and if $K$ is nonempty, 
then (after a rotation) $K$ is the graph of a Lipschitz function $f: \RR\to \RR^2$.  Meeks then showed that
(in that case) $K$ must in fact be a straight line.
\end{comment}

Meeks and Weber~\cite{MeeksWeber}
 showed that every $C^{1,1}$ curve arises as such a blow-up set $K$.
Hoffman and White~\cite{HoffmanWhiteSingularities}
showed that every closed subset of a line arises as such a blow-up set.  
(Kleene~\cite{Kleene} gave another proof of the Hoffman-White result.  Special cases
had been proved earlier by Colding-Minicozzi~\cite{ColdingMinicozziExample}, 
Brian Dean~\cite{BrianDean}, and Siddique Kahn~\cite{SiddiqueKahn}.)

The following questions remain open:
\begin{enumerate}
\item Can $C^1$ in Theorem~\ref{MainFromIntro} be replaced by $C^{1,1}$? The Meeks-Weber examples
show that one cannot prove more regularity than $C^{1,1}$.
\item If $C^1$ can be replaced by $C^{1,1}$,  
does every closed subset of a $C^{1,1}$ curve arise as the blow-up
set $K$ of some sequence $D_n$?  If $C^1$ cannot be replaced by $C^{1,1}$, 
does every closed subset of a $C^1$
curve arise as such a $K$?
\end{enumerate}

\section{Results}\label{Results}

We begin with some definitions.  For simplicity, we work in $\RR^3$, although the results
generalize easily to arbitrary smooth Riemannian $3$-manifolds; see the remark at the 
end of the paper.
A {\bf configuration} is a triple $(U,K,L)$ where $U$ is an open ball in $\RR^3$, an open halfspace
in $\RR^3$ or all of $\RR^3$,
where $K$ is a relatively closed subset of $U$, and where $L$ is a minimal lamination of $U\setminus K$.
Here $K$ should be thought of as a singular set: the configurations $(U,K,L)$ we are most
interested in arise as limits of smooth, properly embedded minimal surfaces, in which case $K$
will be the set of points where the curvature blows up.

We define the {\bf curvature} of a configuration $(U,K, L)$ at a point $p\in L$ to be the norm of the second
fundamental form at $p$ of the leaf that contains $p$.  We define the curvature of the 
configuration $(U,K,L)$ to be $\infty$ at each point of $K$.

A plane $P$ (i.e., a two-dimensional linear subspace of $\RR^3$) is said to be {\bf tangent}
to $(U,K,L)$ at a point $p$ if and only if
\begin{enumerate}
\item $p\in L$ and $P$ is the tangent plane at $p$ to the leaf of the lamination that contains $p$, or
\item $p\in K$.
\end{enumerate}
Thus each point in $L$ has a unique tangent plane, whereas each point in $K$ has (by definition) {\em every} plane
as a tangent plane.

If $(U,K,L)$ is a configuration, the {\bf lift} of $(U,K,L)$ is
\[
\Phi(U,K,L)
=
\{ (x,P): \text{$x\in K\cup L$ and $P$ is tangent plane to $(U,K,L)$ at $x$} \}.
\]
Note that the lift is a relatively closed subset of the Grassmann bundle $U\times G$, where $G$
is the set of all $2$-dimensional linear subspaces of $\RR^3$.
Note also that a configuration is determined by its lift: if $\Phi(U,K,L)=\Phi(U,K',L')$ then $K=K$ and $L=L'$.

\begin{theorem}\label{ConvergenceTheorem}
Let $(U_n, K_n, L_n)$ be a sequence of configurations 
such that $U_n$ converges to a nonempty
open set $U$.  Suppose also that the lifts $\Phi(U_n, K_n, L_n)$ converge
in the Gromov-Hausdorff sense to a relatively closed subset $V$ of $U\times G$.
Then $V$ is the lift $\Phi(U,K,L)$ of a configuration $(U,K,L)$. Furthermore,
\begin{enumerate}
\item For each point $q\in K$, the curvatures of the $(U_n,K_n,L_n)$ blow up at $q$,
meaning that there is a sequence $q_n\in K_n\cup L_n$ such that $q_n$ converges to $q$
and such that the curvature of $(U_n,K_n,L_n)$ at $q_n$ tends to $\infty$ as $n\to\infty$.
\item For each compact subset $C$ of $U\setminus K$, the curvatures of the $(U_n, K_n,L_n)$ are
  uniformly bounded on $C$ as $n\to\infty$.
\item The laminations $L_n$ converge to the lamination $L$ on compact subsets of $U$.
\end{enumerate}
\end{theorem}

Here (and throughout the paper) convergence of open sets $U_n$ to open set $U$
means convergence of  $\RR^3\setminus U_n$
to $\RR^3\setminus U$ in the Gromov-Hausdorff topology.  In particular, if $U_n$ and $U$
are balls, convergence of $U_n$ to $U$ means that the centers and radii of the $U_n$ converge
to the center and radius of $U$.

\begin{proof}
Let $K$ be the set of points $q$ in $U$ such that
\[
   \{q\}\times G \subset V.
\]
First we prove that (1) holds.  For suppose it fails at a point $q\in K$.
By passing to a subsequence, we may assume (for some ball $W$ centered
at $q$) that the curvatures of the $(U_n, K_n, L_n)$ are uniformly bounded
on $W$.  In other words, $W$ is disjoint from each $K_n$ and the curvatures
of the lamination $L_n\cap W$ are uniformly bounded.  By replacing $W$ by
a smaller ball, we can then ensure that the tangent planes to $L_n$
at any two points of $L_n\cap W$ make an angle of at most $\pi/20$ (for example)
with each other.  It follows that if $(x,P)$ and $(x', P')$ are points of $V$
with $x,x'\in W$, then the angle between $P$ and $P'$ is at most $\pi/20$.
But this contradicts the fact that $\{q\}\times G\subset V$, thus proving (1).

Next we prove that (2) holds.  Suppose that $q\in U\setminus K$.  Then
\[
    \{P\in G: (q, P)\in V\} \tag{*}
\]
is a closed subset of $G$ but is not equal to $G$.
Thus there is a closed set $\Sigma\subset G$ with nonempty interior such that
$\Sigma$ is disjoint from the set \thetag{*}.
In other words,
\[
   (\{q\}\times \Sigma) \cap V = \emptyset.
\]
By the Gromov-Hausdorff convergence $\Phi(U_n,K_n, L_n)\to V$, it follows that
there is an open ball $W$ centered at $q$ and compactly contained in $U$ such that
\[
  (\overline{W}\times \Sigma)\cap \Phi(U_n,K_n,L_n) = \emptyset
\]
for all sufficiently large $n$, say $n\ge N$.  It follows immediately that
\begin{enumerate}[\upshape (i)]
\item $K_n\cap W=\emptyset$ for $n\ge N$, and
\item the Gauss map of $L_n\cap W$ omits $\Sigma$ for $n\ge N$.
\end{enumerate}
By a theorem of Osserman~\cite{OssermanCurvatureBound},  (i) and (ii) imply that the curvatures of the $L_n$ are uniformly bounded
(for $n\ge N$) on compact subsets of $W$.  This together with (i) implies that the
curvatures of the $(U_n, K_n, L_n)$ are uniformly bounded on compact subsets of $W$.
This proves (2).

It remains only to prove (3).  Note that the curvature bounds in (2) imply that
every subsequence of the $L_n$ has a further subsequence that converges on
compact subsets of $U\setminus K$ to a lamination $L$ of $U\setminus K$.
But clearly $L$ is determined by $V$.\footnote{In fact, $V\cap ((U\setminus K)\times G)$
is the lift of $(U\setminus K, \emptyset, L)$, so the latter may be recovered from the former
using the projection map from $U\times G$ to $U$.}
Thus the limit $L$ is independent of the subsequence, which means that the
original sequence $L_n$ converges to $L$ on compact subsets of $U\setminus K$.
\end{proof}

We say that configurations $(U_n, K_n, L_n)$ {\bf converge} to configuration $(U,K,L)$
provided $U_n$ converges to $U$ and $\Phi(U_n,K_n,L_n)$ converges in the Gromov-Hausdorff topology
to $\Phi(U,K,L)$.  
From Theorem~\ref{ConvergenceTheorem} together with 
compactness of the space of closed sets under Gromov-Hausdorff convergence, we deduce

\begin{corollary}[Compactness of configurations] 
Suppose $(U_n,K_n, L_n)$ is a sequence of configurations such that $U_n$
converges to a nonempty open set $U$.  Then a subsequence of the $(U_n,K_n,L_n)$
converges to a configuration $(U,K,L)$.
\end{corollary}

A {\bf configuration of disks} is a configuration $(U,\emptyset, L)$ in which each leaf of $L$
is a properly embedded minimal disk in $U$.  We let $\Dd$ be the set of 
all configurations of disks.   We let $\overline{\Dd}$ be the set of all configurations
that are limits of configurations of disks.
Note that $\overline{\Dd}$ is closed under sequential convergence.

\begin{theorem}\label{MainTheorem} Suppose that $(U,K,L) \in \overline{\Dd}$.  
Then $U$ is covered by open balls $\BB$ with the following
properties:
\begin{enumerate}
\item For each point $p\in K\cap \BB$, there is a leaf $L_p$ of $L\cap \BB$ such
that $L_p\cup\{p\}$ is a minimal graph over a planar region and is properly embedded in $\BB$.
\item If $q_n\in K\cap\BB$ converges to $q\in K\cap \BB$, then $L_{q_n}\cup\{q_n\}$ converges
smoothly to $L_q\cup\{q\}$.
\item The singular set $K\cap B$ is contained in a $C^1$ embedded curve $\Gamma$ such that at each 
point $q$ of $K\cap B$, the curve $\Gamma$ is orthogonal to $L_q\cup\{q\}$ at $q$.
\end{enumerate}
\end{theorem}

(See Remark~\ref{Ambient} for the generalization to arbitrary Riemannian $3$-manifolds.)

\begin{proof}
Assertion (1) is due to Colding and Minicozzi~\cite{ColdingMinicozziII}*{Theorem~5.8}.
Assertion (2) follows immediately from Assertion (1).  
To prove Assertion (3), we use the following theorem due to Colding-Minicozzi and Meeks:

\begin{theorem}\label{MeeksTheorem}
If $(\RR^3, K, L)\in \overline{\Dd}$ and if $K$ is nonempty, then $K$ is a line
and the lamination $L$ is the foliation consisting of all planes perpendicular to $L$.
\end{theorem}

(According to \cite{ColdingMinicozziIV}*{Theorem 0.1}, $L$ is a foliation consisting
of parallel planes and $K$ is a Lipschitz curve transverse to those planes.  
According to~\cite{MeeksRegularity}, 
the Lipschitz curve must be a straight line perpendicular to those planes.)

We also use the following proposition, which is a restatement of the $C^1$ case of 
Whitney's Extension Theorem\cite{Whitney}*{Theorem~I}:

\begin{proposition}\label{WhitneyTheorem}
Let $K$ be a relatively closed subset of an open subset $\BB$ of $\RR^n$.
Suppose $\mathcal{V}$ is a continuous line field on $K$, i.e., a continuous function
that assigns to each $p\in K$ a line $\mathcal{V}(p)$ in $\RR^n$.  Suppose also that
if $p_i,q_i\in K$ with $p_i\ne q_i$ converge to $p\in K$, then $\overleftrightarrow{p_iq_i}$
converges to $\mathcal{V}(p)$.   

Then each point $p\in K$ has a neighborhood $W$ such that $K\cap W$ is contained
in the graph $\Gamma$ of a $C^1$ function 
from $\mathcal{V}(p)$ to $(\mathcal{V}(p))^\perp$ such that at each point $q\in W\cap K$, $\mathcal{V}(q)$ is
tangent to $\Gamma$ at $q$.
\end{proposition}

We will apply Proposition~\ref{WhitneyTheorem} with $\mathcal{V}(p) = (\Tan_pL_p)^\perp$.
By Assertion (2) of Theorem~\ref{MainTheorem}, $\mathcal{V}(p)$ depends continuously on $p\in K$.
Let $p_j, q_j\in K\cap \BB$ with $p_j\ne q_j$ converge to $q\in K\cap \BB$. 
 It suffices
to prove that $\overleftrightarrow{p_jq_j}$ converges to $\mathcal{V}(q)$.  

Let $\phi_n:\RR^3\to \RR^3$ be translation by $-q_n$ followed by
dilation by $1/|p_n-q_n|$:
\[
  \phi_n (x) = \frac{x-q_n}{|p_n-q_n|}.
\]
By passing to a subsequence, we may assume that $\phi_n(p_n)$ converges
to a point $p^*$ with $|p^*|=1$.  Thus
\[
   \overleftrightarrow{\mathstrut p_nq_n} = \overleftrightarrow{\mathstrut \phi_n(p_n) \phi_n(q_n)}
   =\overleftrightarrow{\mathstrut \phi_n(p_n) 0} \to \overleftrightarrow{\mathstrut p^*0}.
\]
 Thus it suffices to prove that $\overleftrightarrow{p^*0}$ is equal to $\mathcal{V}(q)$.

Note that $\phi_n(U_n)\to \RR^3$.  Now consider the configurations $(\phi_n(U), \phi_n(K), \phi_n(L))$.
By passing to a further subsequence, we may assume
that these configurations converge  to
a configuration $(\RR^3, K', L')\in \overline{\Dd}$.  Note that $K'$ is nonempty since $0$ and $p^*$
are in $K'$. 
Thus by Theorem~\ref{MeeksTheorem}, $K'$ is a line and $L'$ consists of all planes perpendicular to $K'$.
Since $K'$ contains $0$ and $p^*$, in fact $K'$ is the line through $0$ and $p^*$.

Now by Assertion (2) of the theorem, 
the leaves $\phi_n(L_{q_n}\cup\{q_n\})$ converge smoothly to $\Tan_qL_q$.
Thus $\Tan_qL_q$ is one of the leaves of $L'$, which means that $\Tan_qL_q$ is perpendicular
to $K'$.  In other words, $K'$ is the line $\mathcal{V}(q)$.
\end{proof}

\begin{remark}\label{Ambient}
The definitions and theorems in this paper generalize to arbitrary smooth Riemannian $3$-manifolds.
In particular, Theorem~\ref{MainTheorem} remains true if $U$ is an open geodesic ball of radius $r$ in a $3$-dimensional Riemannian
manifold, provided all the geodesic balls of radius $\le r$ centered at points in $U$ are mean convex.
(This guarantees that if $D$ is a minimal disk properly embedded in $U$, then 
the intersection of $D$ with any geodesic ball in $U$ is a union of disks.)
The proof is almost identical to the proof in the Euclidean case.
\end{remark}

\end{document}